\documentclass{gtart}

\def\ifplaintex{\expandafter\ifx\csname documentclass\endcsname\relax}

\expandafter\ifx\csname epsfbox\endcsname\relax\input epsf\fi

\ifplaintex 
\else
\fi

\def\gt{{\mathsurround=0pt\it $\cal G\mskip-2mu$eometry \&\ 
$\cal T\!\!$opology}}        

\def\gtp{{\mathsurround=0pt\it $\cal G\mskip-2mu$eometry \&\ 
$\cal T\!\!$opology $\cal P\!$ublications}}  

\def\lognumber#1{\def\thelognumber{#1}}
\def\volumenumber#1{\def\thevolumenumber{#1}}
\def\papernumber#1{\def\thepapernumber{#1}}
\def\volumeyear#1{\def\thevolumeyear{#1}}

\def\pagenumbers#1#2{\def\startpage{#1}\def\finishpage{#2}}
\def\published#1{\def\publishdate{#1}}
\def\proposed#1{\def\theproposer{#1}}
\def\seconded#1{\def\theseconders{#1}}
\def\received#1{\def\receiveddate{#1}}

\def\accepted#1{\def\accepteddate{#1}}
\def\asciititle#1{\def\theasciititle{#1}}

\def\asciiaddress#1{\def\theasciiaddress{#1}}

\long\def\asciiabstract#1{\long\def\theasciiabstract{#1}}
\def\asciikeywords#1{\def\theasciikeywords{#1}}

\let\\\par\let\thevolumenumber\relax\let\thepapernumber\relax
\let\thevolumeyear\relax\let\thesamplenumber\relax\let\startpage\relax
\let\finishpage\relax\let\publishdate\relax\let\receiveddate\relax
\let\reviseddate\relax\let\accepteddate\relax\let\theasciititle\relax
\let\theasciiauthors\relax\let\theasciiaddress\relax
\let\theasciiabstract\relax\let\theasciikeywords\relax
\let\theasciiemail\relax\let\theshortauthors\relax\let\theshorttitle\relax

\long\def\maketitlep{   

\count0=\startpage

\gt\hfill      
\hbox to 77pt{\vbox to 0pt{\vglue -15pt\epsfbox{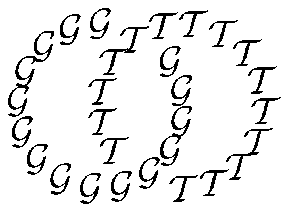}\vss}\hss}
\break
{\small\ifx\thesamplenumber\relax 
Volume \else Sample
\fi\thevolumenumber\ (\thevolumeyear)
\startpage--\finishpage\nl
Published: \publishdate}
\vglue 0.5truein plus 0.4fil minus 0.1truein

{\parskip=0pt\leftskip 0pt plus 1fil\def\\{\par\smallskip}{\ifplaintex\large
\else\Large\fi\bf\thetitle}\par\medskip}   

\vglue 0pt plus 0.1fil 

{\parskip=0pt\leftskip 0pt plus 1fil\def\\{\par}{\sc\theauthors}
\par\medskip}

\vglue 0pt plus 0.1fil 

{\small\parskip=0pt\let\newline\\
{\leftskip 0pt plus 1fil\def\\{\par}{\sl\theaddress}\par}
\expandafter\ifx\theemail\relax    
\relax\else\vglue 5pt plus 0.02fil minus 2pt\def\\{\stdspace{\rm 
and}\stdspace} 
\cl{Email:\stdspace\tt\theemail}\fi
\ifx\theurl\relax                  
\relax\else\vglue 5pt plus 0.02fil minus 2pt\def\\{\stdspace{\rm 
and}\stdspace}
\cl{URL:\stdspace\tt\theurl}\fi\par}

\vglue 7pt plus 0.3fil minus 3pt

{\bf Abstract}
\vglue 5pt plus 0.1fil minus 2pt

\theabstract

\vglue 7pt plus 0.3fil minus 3pt

{\bf AMS Classification numbers}\quad Primary:\quad \theprimaryclass

Secondary:\quad \thesecondaryclass

\vglue 5pt plus 0.3fil minus 2pt

{\bf Keywords:}\quad \thekeywords

\vglue 10pt plus 0.5fil minus 5pt

{\small  Proposed: \theproposer\hfill Received: \receiveddate\nl
Seconded: \theseconders\hfill 
\ifx\reviseddate\relax                         
Accepted: \accepteddate                        
\else
Revised: \reviseddate                          
\fi}
\eject
}       

\let\maketitlepage\maketitlep
\let\maketitle\maketitlepage

\font\phead=cmsl9 scaled 950
\font=cmsl9 scaled 1050
\font\pnum=cmbx10 scaled 913
\font\lnum=cmbx10 
\font\pfoot=cmsl9 scaled 950
\font=cmsl9 scaled 1050
\ifplaintex
\headline{\vbox to 0pt{\vskip -4.5mm\line{\small\phead\ifnum
\count0=\startpage ISSN 1364-0380 (on line)
1465-3060 (printed) \hfill {\pnum\folio}\else\ifodd\count0\def\\{ }%
\ifx\theshorttitle\relax\thetitle\else\theshorttitle\fi\hfill{\pnum\folio}
\else\def\\{ and }{\pnum\folio}\hfill\ifx\theshortauthors\relax\theauthors
\else\theshortauthors\fi\fi\fi}\vss}}
\footline{\vbox to 0pt{\vglue 0mm\line{\small\pfoot\ifnum\count0=\startpage
\copyright\ \gtp\hfill\else
\gt, Volume \thevolumenumber\ (\thevolumeyear)\hfill\fi}\vss
}}
\else
\makeatletter
\def\@oddhead{{\small\lhead\ifnum\count0=\startpage ISSN 1364-0380 (on line)
1465-3060 (printed) \hfill {\lnum\number\count0}\else\ifodd\count0
\def\\{ }\ifx\theshorttitle\relax \thetitle \else\theshorttitle\fi\hfill
{\lnum\number\count0}\else\def\\{ and }{\lnum\number\count0}
\hfill\ifx\theshortauthors\relax 
\theauthors\else\theshortauthors\fi\fi\fi}}\def\@evenhead{@oddhead}
\def\@oddfoot{\small\lfoot\ifnum\count0=\startpage\copyright\ \gtp\hfill\else
\gt, Volume \thevolumenumber\ (\thevolumeyear)\hfill\fi}
\def\@evenfoot{@oddfoot}
\makeatother
\fi

\newwrite\gtoutfile
\long\gdef\makeheadfile{  
{\def\\{, }\def\s{ }
\immediate\openout\gtoutfile head.xxx
\immediate\write\gtoutfile{To: math@arxiv.org}
\immediate\write\gtoutfile{Subject: put OR rep NNNNN:pppp}
\immediate\write\gtoutfile{--text follows this line--}
\immediate\write\gtoutfile{Proxy-for: \ifx\theasciiauthors\relax
\theauthors\else\theasciiauthors\fi\s<\ifx\theasciiemail\relax\theemail\else\theasciiemail\fi>}
\immediate\write\gtoutfile{\noexpand\\}
\immediate\write\gtoutfile{Authors: \ifx\theasciiauthors\relax
\theauthors\else\theasciiauthors\fi}
{\def\\{ }\immediate\write\gtoutfile{Title: \ifx\theasciititle\relax
\thetitle\else\theasciititle\fi}}
\immediate\write\gtoutfile{Subj-class: GT or GR or SG or ...}
\immediate\write\gtoutfile{MSC-class: \theprimaryclass\ifx\thesecondaryclass\relax\else, \thesecondaryclass\fi}
\immediate\write\gtoutfile{Journal-ref: Geom. Topol. \thevolumenumber\s
(\thevolumeyear) \startpage-\finishpage}
\immediate\write\gtoutfile{Comments: Published in Geometry and Topology at}
\immediate\write\gtoutfile{    http://www.maths.warwick.ac.uk/gt/GTVol\thevolumenumber/paper\thepapernumber.abs.html}
\immediate\write\gtoutfile{\noexpand\\}
\immediate\write\gtoutfile{}
\ifx\theasciiabstract\relax
\immediate\write\gtoutfile{\theabstract}\else
\immediate\write\gtoutfile{\theasciiabstract}\fi
\immediate\write\gtoutfile{}
\immediate\write\gtoutfile{\noexpand\\}
\immediate\write\gtoutfile{}
\immediate\closeout\gtoutfile}}  

\def\maketitlepage{\maketitlep\makeheadfile}
\let\maketitle\maketitlepage

\lognumber{313}
\volumenumber{7}\papernumber{9}
\volumeyear{2003}
\pagenumbers{311}{319}
\received{28 March 2003}
\accepted{30 April 2003}
\published{17 May 2003}
\proposed{Stephen Ferry}
\seconded{Ronald Stern, Benson Farb}

\usepackage{amsmath,amssymb}
\input xypic

\def\bignohang#1{\item}
\def\lspec{{\mathbb L}}
\def\headroom#1{\raise #1ex\hbox{\vphantom{.}}}   
\def\slnr{\hbox{SL}_n(\reals)}

\def\trace{{\hbox{tr}\,}}
\def\ker{{\hbox{ker}\,}}
\def\coker{{\hbox{coker}\,}}
\def\ind{{\hbox{ind}\,}}
\def\dslash{\,\hbox{\slash}\kern-8.0pt D}
\def\tildedslash{\widetilde{\dslash}}

\def\tinysig{\hbox{\tiny sig}}
\def\tinysigtwo{{\tiny\hbox{\sig}_{\!(2)\headroom{2.0}}}}
\def\integers{{\mathbb{Z}}}

\def\naturals{{\mathbb N}}
\def\rationals{{\mathbb Q}}
\def\reals{{\mathbb R}}
\def\realsp{\reals{\mathbb P}}
\def\sig{{\hbox{sig}\,}}
\let\bigskip\medskip

\begin{document}

\title{On Invariants of Hirzebruch and Cheeger--Gromov}
\asciititle{On Invariants of Hirzebruch and Cheeger-Gromov}

\authors{Stanley Chang\\Shmuel Weinberger}

\address{Department of Mathematics, Wellesley College\\Wellesley, MA
02481, USA\\\smallskip\\{\rm and}\\\smallskip\\Department 
of Mathematics, University of
Chicago\\Chicago, IL 60637, USA}

\asciiaddress{Department of Mathematics, Wellesley College\\Wellesley, MA
02481, USA\\and\\Department of Mathematics, University of
Chicago\\Chicago, IL 60637, USA}

\email{shmuel@math.uchicago.edu, sschang@palmer.wellesley.edu}

\asciiabstract{We prove that, if M is a compact oriented manifold of
dimension 4k+3, where k>0, such that pi_1(M) is not torsion-free, then
there are infinitely many manifolds that are homotopic equivalent to M
but not homeomorphic to it.  To show the infinite size of the
structure set of M, we construct a secondary invariant tau_(2):
S(M)-->R that coincides with the rho-invariant of Cheeger-Gromov. In
particular, our result shows that the rho-invariant is not a homotopy
invariant for the manifolds in question.  }

\begin{abstract}
We prove that, if $M$ is a compact oriented manifold
of dimension $4k+3$, where $k>0$, such that $\pi_1(M)$
is not torsion-free, then there are infinitely many manifolds
that are homotopic equivalent to $M$ but not homeomorphic to it.
To show the infinite size of the structure set of $M$, 
we construct a secondary invariant $\tau_{(2)}\colon S(M)\to \reals$
that coincides with the $\rho$--invariant of Cheeger--Gromov. In
particular, our result shows that the $\rho$--invariant is not
a homotopy invariant for the manifolds in question.
\end{abstract}

\asciikeywords{Signature, L^2-signature, structure set, rho-invariant}


\maketitlepage


In their analysis of free involutions on the sphere, Browder and
Livesay \cite{BL} introduced an invariant which they used to study the
existence and uniqueness question for desuspensions. In particular,
they produced infinitely many manifolds homotopy equivalent, but not
diffeomorphic or even homeomorphic, to $\realsp^{4k-1}$ for $k>1$.
Hirzebruch \cite{Hi} soon gave an alternative definition of this invariant
as follows: Any manifold $M$ of dimension $4k-1$ with
$\pi_1(M)=\integers_2$ has a multiple $\ell M$ which bounds some $W^{4k}$
such that $\pi_1(M)\to \pi_1(W)$ is an isomorphism. Define
$$\tau(M)=\frac{1}{\ell}\cdot g-\sig(\widetilde{W})
   =\frac{1}{\ell}\,\left(\,\trace g\,\vert_{H_+}-\trace g\,\vert_{H_-}\right)=
\frac{1}{\ell}\,\left(\,\sig(\widetilde{W})-2\,\,\sig(W)\right),$$
where $g\in\integers_2$ is the nontrivial element and $\widetilde{W}$
is the two-fold cover of $W$ with $H_\pm\equiv H_\pm^{2k}(\widetilde{W})$.
The fact that $\tau(M)$ is independent of $W$ can be seen either using
the $G$--signature formula \cite{APS} or by a bordism argument.

\bigskip\noindent
{\bf Remark}\qua Actually, Hirzebruch's invariant was only shown a bit
later to be equal (up to a sign) to the Browder--Livesay invariant
by Hirzebruch and Lopez de Medrano \cite{LdM}.

\bigskip
It is quite simple to modify $\tau$ to obtain invariants of arbitrary
$M^{4k-1}$ with finite fundamental group $G=\pi_1(M)$
by using the $G$--signature of $\widetilde{W}$
modulo the regular representation (as in \cite{Wa}) or by using characters
equal to $\trace g\,\vert_{H_+}-\trace g\,\vert_{H_-}$ for any $g\ne e$.
As a very special and important case, one can
even simply use $\sig(\widetilde{W})-|G|\,\sig(W)$ which equals
$\sum_{g\ne e}\chi_g(\sig_G(W))$ to obtain a quite useful invariant.

As an application of their index theorem for manifolds with boundary,
Atiyah, Patodi and Singer \cite{APS} obtain such invariants of $M$ without
introducing $W$. More precisely, the APS invariant requires the use 
of an auxiliary Riemannian metric on $M$, but their index theorem for 
manifolds with boundary shows that the number obtained is 
independent of this metric.  All of the numerical invariants 
mentioned above are differences of $\eta$--invariants (which measure 
``spectral asymmetry'') for the signature operator twisted by 
suitable flat bundles.

The APS method admits extensive generalization to the case 
in which $\pi_1(M)$ is
possibly infinite.  In \cite{A} Atiyah generalized the index theorem from compact 
closed manifolds to universal covers of compact manifolds (or even manifolds with 
proper actions of a discrete group). Using von Neumann traces, he
defined $$\ind_\Gamma(\tildedslash)=\dim_\Gamma(\ker\tildedslash\hbox{~on~}
L^2(\widetilde{M})) - \dim_\Gamma(\coker\tildedslash\hbox{~on~}
L^2(\widetilde{M})).$$  

Both \cite{L1} and \cite{L2} are  excellent references for ideas related to the von 
Neumann trace and its applications in geometry and topology.

\bigskip\noindent
{\bf Note}\qua 
Atiyah showed that 
$\ind_\Gamma(\tildedslash)
=\ind(\dslash)$, which generalizes the {\it multiplicativity}
property of ordinary indices for finite coverings (see appendix). 
This result is quite remarkable since, if the group is infinite,
the elements of $\ker\dslash$ never lift to be $L^2$ on the universal cover,
but they nonetheless lead indirectly to
elements of the $L^2$ kernel (at the level of indices, of course).
Cheeger and Gromov
\cite{CG1,CG2} introduced the analogue of the APS invariant in this
setting: namely they studied an $L^2$--invariant $\eta_{(2)}(\widetilde{M})$
of compact $M$ 
and the difference $\rho_{(2)}(M)=\eta_{(2)}(\widetilde{M})-\eta(M)$ which
is metric independent. For $\Gamma=\integers_2$ the quantity $\rho_{(2)}(M)$
is simply the Hirzebruch invariant or the Browder--Livesay invariant.

 Thom's classical work on cobordism implies 
that every compact odd dimensional oriented manifold $M$
has a multiple $rM$ which is the boundary of an oriented manifold $W$. 
Hausmann \cite{Ha} showed furthermore that, for every such $M$ with fundamental group
$\Gamma$, there is a manifold $W$ such that 
$\Gamma$ injects into $\pi_1(W)$ and
$\partial W=rM$, for some multiple $rM$ of $M$; that is, he showed that 
nullcobordant manifolds bound in such a way that their fundamental 
groups inject. 
If $\Gamma$ is a group, we
can define the $L^2$--signature $\sig_{(2)}^\Gamma
(V)=\dim_\Gamma(V^+)
-\dim_\Gamma(V^-)$ for any $\ell^2(\Gamma)$--module $V$
endowed with nonsingular symmetric bilinear form. 
As usual, one uses an inner product to
rewrite the symmetric bilinear form in terms of a self-adjoint
operator, and then uses the spectral theorem to obtain projections to the
positive and negative definite part. If $N^{4k}$ is a $\Gamma$--space, let
$\sig_{(2)}^\Gamma (N)$  be given by the $L^2$--signature of the symmetric
form induced by cap product on its middle cohomology.

 Define a new ``Hirzebruch type'' invariant $\tau_{(2)}$
given by $$\tau^G_{(2)}(M) = \frac{1}{r}\left(\sig^G_{(2)}(W_G)-\sig(W)
\right),$$
where $G=\pi_1(W)$ and $W_G$ is the induced $G$--cover of $W$.  
Suppose that $rM=\partial W$ with $\pi_1(M)=\Gamma$
injecting into $\pi_1(W)=G$. If $G$ injects into some larger group
$G'$, and let $W_{G'}$ be the $G'$--space induced from the $G$ action on 
$W_G$ to $G'$.  Notice that $\partial(W_{G'}/G')=rM$, so we can use the 
larger group $G'$ to define  $\tau_{(2)}$.  However, 
by the $\Gamma$--induction property of Cheeger--Gromov
\cite[page 8, equation (2.3)]{CG2}, we have
$$
\begin{array}{rcl}
\tau_{(2)}^{G'}(M)&=& \frac{1}{r}\big(
    \sig_{(2)}^{G'}(W_{G'})-\sig(W_{G'}/G')\big)\\
  &=& \frac{1}{r}\big(\sig_{(2)}^{G\headroom{2}}(W)-\sig(W)\big)\\
  &=& \tau_{(2)}^{G\headroom{2}}(M).
\end{array}
$$
So one can pass to any larger group without changing the value of the 
invariant. Now, given two manifolds $W$ 
and $W'$ with the required bounding properties,
we can use the large group $G'=\pi_1(W)\ast_\Gamma\pi_1(W')$, which 
contains both fundamental groups, and the usual Novikov 
additivity argument, to see that $\tau_{(2)}$ is independent of all 
choices.

In fact, according to \cite{LS1}, we have $\tau_{(2)}=\rho_{(2)}$, but
we do not need this equality for the proofs of our main theorems.
The goal of this paper is to prove a few theorems about
$\tau_{(2)}(M)$ (or
$\rho_{(2)}(M)$) that extend or are analogues of the classically
understood situation
for finite fundamental group.

\bigskip\noindent
{\bf Theorem 1}\qua\sl
For any compact oriented manifold $M^{4k-1}$,
where $k\ge2$,  such
that $\pi_1(M)$ is not torsion-free,
there exist 
infinitely many manifolds
$M_i$ which are 
simple homotopy equivalent and tangentially equivalent
to $M$, but not homeomorphic to $M$. If $M$ is smooth,
the $M_i$ can be taken to be smooth as well. 
This infinite number is detected by $\tau_{(2)}$.\rm

\bigskip\noindent
{\bf Remarks}\qua

\begin{enumerate}
\bignohang
{1.} Very similar arguments to the proof of Theorem 1
show that any $M^{4k+2}$
with non-torsion-free fundamental group has infinite diffeotopy group.

\bignohang
{2.} The analogue for $M^{4k+1}$ is false. For instance, the real projective
space $\realsp^{4k+1}$ has a finite structure set. See \cite{LdM}.
Also, we cannot dispense with the
orientability hypothesis as $\realsp^2 \times S^5$ has finite structure set.

\bignohang
{3.} There are predictions in the literature about $L_{4k}(\Gamma)\otimes
\rationals$ and its expression in terms of the group homology of normalizers
of finite subgroups. See \cite{FJ} and page 261 of \cite{We2}. 
These would suggest that Theorem 1 is
true except for the thorny problem of self-homotopy equivalences of $M$. These
maps $h:M\to M$ could represent nontrivial elements of the structure
set $S(M)$, but nonetheless they are {\it not} the elements described
in Theorem 1. 

\bignohang
{4.} If $M$ is smooth, all the $M_i$ can be taken smooth as well, without
even being homeomorphic. The extension of $\tau_{(2)}$ from smooth to
topological manifolds is not very hard by bordism methods. 
After multiplication by a positive integer, every
topological manifold is cobordant to a smooth manifold $N$ 
via a cobordism $V$ (all with the same
fundamental group), so we can consider $\tau_{(2)}(N)  -
(\sig_{(2)}-\sig)(V)$.

\bignohang
{5.} The theorem above asserts that, whenever $\pi_1(M)$ has torsion,
then $\tau_{(2)}$ is not a homotopy invariant. Mathai \cite{Ma} conjectured that
the converse is true for $\rho_{(2)}$
and proved it for torsion-free crystallographic 
fundamental group. Thus, the above result is the converse to Mathai's
conjecture.
\end{enumerate}

Mathai's conjecture is an analogue of a conjecture and theorem of the second
author \cite{We1} for Atiyah--Patodi--Singer invariants. Keswani \cite{K} proved
a general homotopy invariance theorem:
if the assembly map $K_\ast(B\Gamma)\to K_\ast(
C^\ast_{\hbox{\tiny{max}}\headroom{-1}}\Gamma)$ is surjective, 
then $\rho_{(2)}$ is a homotopy
invariant for manifolds with fundamental group $\Gamma$.
This hypothesis holds for fundamental group of real and complex hyperbolic
manifolds for torsion-free amenable groups by the theorem of Higson
and Kasparov \cite{HK}. However, it fails for all nontrivial groups satisfying
property $T$, ie, for many torsion-free lattices in higher rank
or in $\hbox{Sp}\,(n,1)$. Our next theorem somewhat
repairs this problem.

\bigskip\noindent
{\bf Theorem 2}\qua\sl If $\Gamma$ is a torsion-free discrete subgroup
of $\slnr$, then $\rho$ is a homotopy invariant for manifolds with 
fundamental group \( \Gamma \). More generally, if $\Gamma=\pi_1(M)$ is 
residually
finite and the Borel conjecture holds for $\Gamma$, then $\rho_{(2)}$ is
a homotopy invariant (see eg \cite{FJ}).\rm

\bigskip
In a future paper, the first author plans to remove the residual finiteness
hypotheses arising in this work.

We now turn to the proof of Theorem 1. Let $M$ be as in the theorem
and $\integers_n\subset\pi_1(M)=\Gamma$ be a nontrivial cyclic subgroup. The
$M_i$ are just the results of acting by elements of $\hbox{im}\,
({L}_{4k}
(\integers_n)\to {L}_{4k}(\Gamma))$
on $S(M)$ in the surgery exact sequence
$$\cdots\to {H}_{4k}(M,\lspec)\to {L}_{4k}(\Gamma)
\to S(M)\to {H}_{4k-1}(M,\lspec)\to {L}_{4k-1}(
\Gamma)\to
\cdots$$

The image of $A: {H}_{4k}(M,\lspec)\to {L}_{4k}
(\Gamma)$ consists of quadratic forms arising from closed manifolds
and for such manifolds we have $\sig_{(2)}^\Gamma(N)=\sig(N)$ by Atiyah \cite{A}. 
Note that this higher signature may be viewed as a map 
$\sig_{(2)}^{\Gamma\headroom{1.5}}:
L_{4k}(\Gamma)\to\reals$ (see \cite{Wa} for the
definition of $L$--groups). This
image must equal $\hbox{ker}\,({L}^{\headroom{1.5}}_{4k}
  (\Gamma)\to S(M))$,
the set of all $a\in {L}_{4k}(\Gamma)$ for which $a(M)=M$
diffeomorphically.\footnote{More precisely, the equality 
represents a diffeomorphism
homotopic to the homotopy equivalence implicit in the fact
that $a(M)$ is a member
of a structure set.}
 If $a\in {L}^{\headroom{1.5}}_{4k
\headroom{-.7}}
  (\Gamma)$, let
$V_a$ be the form representing it. Let $\alpha_\Gamma:
{L}^{\headroom{1.5}}_{4k}(\Gamma)\to\reals$ be given by
the homomorphism
$\alpha_\Gamma(V)\equiv \sig_{(2)\headroom{-1}}^\Gamma(V)-\sig(V)$. 
Hence, if $a(M)=M$, then $\alpha_\Gamma(V_a)=0$. We would like to
produce infinitely many
$V_{a'}$ such that $\alpha_\Gamma(V_{a'})\ne0$.

To show that the homomorphism $\alpha_\Gamma:
{L}_{4k}(\Gamma)\to\reals$ has infinite image,
it suffices to show that it is nontrivial.
Let $p_\ast: {L}_{4k}(\integers_n)\to {L}_{4k}
(\Gamma)$ be the map induced by the injection $p:\integers_n\to\Gamma$.
By the induction property of $\Gamma$--dimension, we have
the commutative diagram
$$
\xymatrix@R=8ex@C=4em{
{L}_{4k}(\integers_n)
     \ar[d]_{p_\ast}\ar[r]^-{\tinysigtwo^{\hskip -.15in 
\integers_{n\headroom{-.5}}}}&\reals\\
{L}_{4k}(\Gamma)\ar[ur]_-{\tinysigtwo^{\hskip -.15in
\Gamma_{\headroom{-.4}}}}&
}
$$
so we only have to check that $\alpha_{\integers_n}: 
{L}_{4k}(\integers_n)\to\reals$ is nontrivial.  
We note that, inverting $2$,
it makes no difference whether the coefficient ring is 
$\integers$ or $\rationals$ \cite{Ran}
or whether one works with projective, free or based-free modules (by 
the well-known ``Rothenberg exact sequences").  Thus, we need only
 produce a bilinear form on a projective 
module over $\rationals\integers_n$, namely
 the one-by-one bilinear form $[1]$ on the 
module $\rationals$, which is projective
over $\rationals\integers_n$.  
By definition, our invariant is $\frac{1}{n} - 1$ for this element, which is
always nontrivial for $n>1$.

If $rM^{\headroom{1.5}}$ 
cobounds some space $W$ for which $\pi_1(M)$ injects into
$G=\pi_1(W)$, we can extend each boundary component $M$ by the cobordism
$Y$, so that $rM'$ cobounds $W'$, where $W'$ is obtained from $W$ and
$r$ copies of $Y$. 
Hence by Novikov additivity and $\Gamma$--induction, we have
$${
\begin{array}{rcl}
  \tau_{(2)}^G(M)-\tau_{(2)}^G(M') & = & 
             \frac{1}{r}\big(\sig_{(2)}^{G\headroom{2}}(W_G)-\sig(W)\big)
            -\frac{1}{r}\big(\sig_{(2)}^G
                  (W'_G)-\sig(W')\big)\\
   &=& \frac{1}{r}\big(\sig_{(2)}^{G\headroom{2}}(W_G)-\sig(W)\big)-
   \frac{1}{r}\big(\sig_{(2)}^G
                  (W))\\
    &&\hspace{2.5cm}+\,r\cdot\sig_{(2)}^G(Y_G-\sig(W)-r\cdot\sig(Y)\big)\\
    &=& \sig_{(2)}^{G\headroom{2}}(Y_G)-\sig(Y)\\
    &=& \sig_{(2)\headroom{-2}}^{\Gamma\headroom{2}}(Y_\Gamma)-\sig(Y)\\
     &=&\alpha_\Gamma(V), 
\end{array}  }
$$
where $V\in {L}_{4k}(\Gamma)$ represents the cobordism
$Y$. 

Let $V_a\in {L}_{4k}(\Gamma)$ such that $\alpha_\Gamma(V_a)
\ne0$, where $\Gamma=\pi_1(M)$. Since $(M\to M)\in S_{4k-1}(M)$,
there is a cobordism $Y_a$ such that 
$\partial^{\headroom{1.5}} Y_a=M\cup M'$ representing
$V_a$. Clearly $M$ and $M'$ are homotopy equivalent but not
diffeomorphic.
Consider $V_{a_i\headroom{-.7}}^{\headroom{1.5}}
\in  {L}_{4k}(\Gamma)$ for all $i\in\naturals$ such
that $\alpha_\Gamma(V_{a_i})\ne\alpha_\Gamma(V_{a_j})$ for all $i\ne j$ and
all nonzero. Let $M_i$ be the result of acting on $M$ by $V_{a_i}$. Then 
$\tau_{(2)}^{G\headroom{1.5}}(M)
-\tau_{(2)}^G(M_i)=\alpha_\Gamma(V_{a_i})\ne0$. Also
$$
\begin{array}{rcl}
  \tau_{(2)}^G(M_i)-\tau_{(2)\headroom{-2}}^G(M_j)
      &=& \tau_{(2)}^G(M_i)-\tau_{(2)}^G
    (M)+\tau_{(2)}^G(M)-\tau_{(2)}^G(M_j)\\ 
                 &=& -\alpha_\Gamma(V_{\alpha_i})+\alpha_\Gamma(V_{\alpha_j})\ne0,
\end{array}
$$
as desired. \endproof

Now we prove Theorem 2. For a finite group 
$\Gamma$, the rho-invariant
$\rho_{(2)} =\eta_{(2)}- \eta$ is equal to
   $\frac{1}{|\Gamma|}\,\eta_{\otimes\reals \Gamma}-\eta$,
where $\eta_{\otimes\reals \Gamma}$ is the eta--invariant for $M$ with
coefficients in the flat bundle associated to the regular representation.
According to \cite{LS2}, if $\Gamma$ is residually finite with descending
quotients $\Gamma_k$, then 
$$\rho_{(2)}=\lim_{k\to\infty}\left(\frac{1}{|\Gamma_k|}
\left(\eta_{\otimes\reals \Gamma}-\eta\right)\right).$$
Indeed, Cheeger and Gromov introduced $\rho_{(2)}$ to generalize the
right-hand side to non-residually finite fundamental group. However,
according to \cite{We2}, for {\it any} representation $V$ of a group
satisfying the Borel conjecture, the quantity $\frac{1}{\dim V}\,
\eta_{\otimes V}-\eta$ is a homotopy invariant. Since $\rho$ is the
limit of homotopy invariants, a fortiori, it too is a homotopy invariant. 
\endproof

If $\Gamma$ is a discrete subgroup of $\slnr$, then it is a classical lemma
of Selberg that $\Gamma$ is residually finite. The fact that the Borel
conjecture holds if $\Gamma$ is in addition torsion-free is the main
result of \cite{FJ}.

\bigskip\noindent
{\bf Problem}\qua If $\pi_1(M)$ has many conjugacy classes of elements of
finite order, then one would expect that $S(M)$ is rather larger than
the $\integers$ we just detected. Can this statement
be proven unconditionally? Wolfgang L\"uck has pointed out that, using
the traces associated to elements of finite order that have finitely
many conjugates (eg, central elements), one can extend Theorem 1 and
detect a larger $S(M)$. But this result is still rather smaller
than the predicted size of the structure set.

\rk{Appendix: Acyclic groups and Atiyah's theorem}

For the convenience of the reader we shall prove here the theorems
of Atiyah and Hausmann invoked in the paper. Let $M$ be a closed manifold.
According to \cite{BDH}, there is an acyclic group $A$ with an injective
homomorphism $h\colon\pi_1(M)\to A$. Since $A$ is acyclic, the map
$\Omega_n(K(A,1))\to\Omega_n(\ast)$ is an isomorphism, so that, if
$M$ bounds (ie, $M$ is zero in $\Omega_n(\ast)$), it bounds a manifold
$W$ such that
$$
\xymatrix{
  M \ar[d]\ar[r]
  &K\bigl(\pi_1(M),1\bigr)\ar[d]\\
  W\ar[r]
  &K(A,1)\\}
$$
commutes; {\it a fortiori\/} $\pi_1(M)\to \pi_1(W)$ is injective. This
argument gives Hausmann's result.

Now for Atiyah's theorem we observe, by the induction and bordism
invariance properties of $\sig_{(2)}$, that we have the commutative diagram
$$
\xymatrix@R=8ex@C=4em{
  & \Omega_{4k}(\ast)\ar[dl]_-\tinysig\ar[d]\ar[dr]^-\tinysigtwo&\\
  \integers
  &\Omega_{4k}\bigl(K(\pi,1)\bigr)\ar[l]^(.6)\tinysig\ar[d]
  \ar[r]_(.6){\tinysigtwo^{\hskip -.3cm \pi}}
  &\reals\\
  &\Omega_{4k}\bigl(K(A,1)\bigr)
  \ar[ul]^-\tinysig\ar[ur]_-{\tinysigtwo^{\hskip -.3cm A}}
}
$$
Since $\Omega_{4k}(\ast)\to \Omega_{4k}(K(A,1))$ is an isomorphism
and $\sig=\sig_{(2)}$ for the trivial group, the conclusion follows.

\rk{Acknowledgements} We would like to thank both Wolfgang L\"uck and
the referee for their helpful comments.

Stanley Chang was partially supported by NSF Grant DMS-9971657;
Shmuel Weinberger was also partially supported by an NSF grant.

\end{document}